\documentclass[12pt,a4paper]{amsart}
\usepackage{amssymb,amsmath}

\headsep=1cm \numberwithin{equation}{section}
\newtheorem{theorem}{Theorem}[section]
\newtheorem{lemma}[theorem]{Lemma}
\newtheorem{proposition}[theorem]{Proposition}
\newtheorem{corollary}[theorem]{Corollary}

\theoremstyle{definition}

\theoremstyle{remark}

\newtheorem{example}[theorem]{Example}

\newcommand{\Hom}{\operatorname{Hom}}

\newcommand{\Ann}{\operatorname{Ann}}

\newcommand{\Max}{\operatorname{Max}}

\newcommand{\BN}{\Bbb N}

\newcommand{\lo}{\longrightarrow}
\newcommand{\fm}{\frak{m}}
\newcommand{\fp}{\frak{p}}

\newcommand{\fa}{\frak{a}}

\begin{document}
\author[Divaani-Aazar and Mafi ]{Kamran Divaani-Aazar and Amir Mafi }
\title[A new characterization of commutative Artinian rings]
{A new characterization of commutative Artinian rings}

\address{K. Divaani-Aazar, Department of Mathematics, Az-Zahra University,
Vanak, Post Code 19834, Tehran, IRAN and Institute for Studies in
Theoretical Physics and Mathematics, P.O.Box 19395-5746, Tehran,
Iran.} \email{kdivaani@ipm.ir}

\address{A. Mafi, Institute of Mathematics, University for Teacher
Education, 599 Taleghani Avenue, Tehran 15614, Iran.}



\thanks{This research was in part supported by a grant from IPM, Tehran IRAN}

\begin{abstract}
Let $R$ be  a commutative Noetherian ring. It is shown that $R$ is
Artinian if and only if every $R$-module is good, if and only if
every $R$-module is representable. As a result, it follows that
every nonzero submodule of any representable $R$-module is
representable if and only if $R$ is Artinian. This provides an
answer to a question which is investigated  in [{\bf 1}].
\end{abstract}

\maketitle

\section{Introduction}

All rings considered in this paper are assumed to be commutative
with identity. There are several characterizations of  Artinian
rings. In particular, it is known that a Noetherian ring $R$ is
Artinian if and only if every prime ideal of $R$ is maximal. In
this article, we present a new characterization of Artinian rings
according to the notions of primary decomposition and (its dual)
secondary representation. To do so, we need to introduce a
generalization of the notions Hopficity and co-Hopficity.

In [{\bf 2}] V. A. Hiremath introduced the concept of Hopficity
for $R$-modules. The dual notion is defined by K. Varadarajan
[{\bf 6}]. An $R$-module $M$ is said to be Hopfian (resp.
co-Hopfian) if any surjective (resp. injective) $R$-homomorphism
is automatically an isomorphism. We refer the reader to [{\bf 6}]
for reviewing the most important properties of Hopfian and
co-Hopfian $R$-modules. We extend these definitions as follows: An
$R$-module $M$ is said to be semi Hopfian (resp. semi co-Hopfian)
if for any $x\in R$, the endomorphism of $M$ induced by
multiplication by $x$ is an isomorphism, provided it is surjective
(resp. injective). Clearly any Hopfian (resp. co-Hopfian)
$R$-module is semi-Hopfian (resp. semi co-Hopfian). Also, it is
obvious that $R$, as an $R$-module, is Hopfian (resp. co-Hopfian)
if and only if it is semi Hopfian (resp. semi co-Hopfian).

As the main result of this note, we establish the following
characterization of Artinian rings.

\begin{theorem} Let $R$ be a commutative Noetherian ring. Then the
following are equivalent:\\
i) $R$ is Artinian.\\
ii) Every nonzero $R$-module is good.\\
iii) Every $R$-module is semi Hopfian.\\
iv) Every nonzero $R$-module is representable.\\
iv') Every nonzero Noetherian $R$-module is representable.\\
v) Every $R$-module is semi co-Hopfian.\\
v') Every Noetherian $R$-module is semi co-Hopfian.\\
v'') Every Noetherian $R$-module is co-Hopfian.\\
vi) Every nonzero $R$-module is Laskerian.

\end{theorem}

In [{\bf 1}] the following question was investigated: When are
submodules of representable $R$-modules representable? In that
paper [{\bf 1}, Theorem 2.3], it is shown that this is the case,
when $R$ is Von Neumann regular. For a Noetherian ring $R$, we
prove that every nonzero submodules of any representable
$R$-modules is representable if and only if $R$ is Artinian (see
2.4).

\section{The proof of the main theorem}

Recall that a nonzero $R$-module $M$ is called {\it good}, if its
zero submodule possesses a primary decomposition. A nonzero
$R$-module $S$ is said to be {\it secondary}, if for any $x\in R$,
the map induced by multiplication by $x$ is either surjective or
nilpotent. We say the $R$-module $M$ is {\it representable}, if
there are secondary submodules $S_1, S_2,\dots ,S_k$ of $M$ such
that $M=S_1+ S_2+\dots +S_k$. The two notions primary
decomposition and secondary representation are dual concepts. We
refer the reader to [{\bf 3}, Appendix to \S6], for more details
about secondary representation. Also, recall that an $R$-module
$M$ is said to be {\it Laskerian}, if any submodule of $M$ is an
intersection of a finite number of primary submodules.

\begin{lemma} i) Every finitely generated $R$-module is
Hopfian. \\
ii) Every Artinian $R$-module is co-Hopfian.\\
iii) Every good $R$-module is semi Hopfian.\\
iv) Every representable $R$-module is semi co-Hopfian.
\end{lemma}

{\bf Proof.} i) See [{\bf 7}, Proposition 1.2].\\
ii) is well known and can be checked easily. \\
iii) Let $x\in R$ be an $M$-coregular element of $R$. Let
$0=\bigcap_{i=1}^n Q_i$ be a primary decomposition of the zero
submodule of $M$. Fix $1\leq i\leq n$. Since $Q_i$ is a proper
submodule of $M$ and $\frac{M}{Q_i}\overset{x}{\lo}\frac{M}{Q_i}$
is either injective or nilpotent, it follows that $x$ is
$\frac{M}{Q_i}$-regular. Now, if $xm=0$ for some element $m$ in
$M$, then for each $i$, it follows that $xm\in Q_i$  and so $m\in
Q_i$. Hence
$m=0$, and so $x$ in $M$-regular as required.\\
iv) is similar to (iii). $\Box$

\begin{example} Let $N$ be a nonzero co-Hopfian $R$-module. Set
$M=\bigoplus_{i\in\BN} N$. Then $M$ is semi co-Hopfian, but it is
not co-Hopfian. To this end define the $R$-homomorphism $\psi:M\lo
M$ by $\psi(m_1,m_2,\dots)=(0,m_1,m_2,\dots)$ for all
$(m_1,m_2,\dots)\in M$. Then $\psi$ is  injective, while it is not
surjective.
\end{example}

{\bf Proof of theorem 1.1.} $(i)\Rightarrow (ii)$ Let $M$ be an
$R$-module. Since $R$ is Artinian, it is representable as an
$R$-module. Hence $M\simeq \Hom_R(R,M)$ is good, by [{\bf 4},
Theorem 2.8]. The implications $(ii)\Rightarrow (iii)$,
$(iv)\Rightarrow (v)$ and $(iv)'\Rightarrow (v)'$ follow, by 2.1.

Now we prove $(iii)\Rightarrow (v)$. Let $M$ be an $R$-module and
$D(\cdot)=\Hom_R(\cdot,E)$, where $E$ is an injective cogenerator
of $R$. Let $x\in R$ be such that the map $M\overset{x}{\lo} M$ is
injective. Then the map $D(M)\overset{x}{\lo}D(M)$ is surjective
and it is also injective, because $D(M)$ is semi Hopfian. But this
implies that $x$ is $M$-coregular, as the functor $D(\cdot)$ is
faithfully exact. Hence $M$ is semi co-Hopfian.

$(v)'\Rightarrow (i)$ Suppose the contrary and assume that
$\fp\subset\fm$ is a strict containment of prime ideals of $R$.
Let $x\in\fm\smallsetminus\fp$. Then $x$ is $R/\fp$-regular, but
it is not $R/\fp$-coregular. We achieved at a contradiction.
Therefore every prime ideal of $R$ is maximal and so $R$ is
Artinian.

Next, we prove $(i)\Rightarrow (iv)$. Since $R$ is Artinian,
$\Max(R)$ is finite. Let $\Max(R)=\{\fm_1,\dots,\fm_k\}$. There
are $\fm_i$-primary ideals $\fa_i$ of $R$ such that
$R\simeq\prod_{i=1}^k R/\fa_i$. Let $F=\bigoplus_{j\in J}R$ be an
arbitrary free $R$-module. Set $S_i=\bigoplus_{j\in J} R/\fa_i$
for $i=1,2,\dots,k$. Then
$$F\simeq\bigoplus_{j\in J} (\prod_{i=1}^k R/\fa_i)\simeq \prod_{i=1}^k
S_i=\bigoplus_{i=1}^k S_i.$$ It follows that for each $i=1,2,\dots
,k$, the $R$-module $S_i$ is $\fm_i$-secondary and hence $F$ is
representable. But any $R$-module is homeomorphic image of some
free $R$-module and so the conclusion follows. Note that one can
check easily that any nonzero quotient of a representable
$R$-module is also representable.

It follows from [{\bf 8}, Theorem] that the statements (i) and
(v'') are equivalent. Let $N$ be a proper submodule of an
$R$-module $M$. Then $N$ possesses a primary decomposition if and
only if the $R$-module $M/N$ is good. Thus (ii) and (vi) are
equivalent. Now, because the implications $(iv)\Rightarrow (iv')$
and $(v)\Rightarrow (v')$ are clearly hold, the proof is complete.
$\Box$

\begin{corollary}  Let $M$ be an $R$-module such that the ring
$R/\Ann_RM$ is Artinian. Then $M$ is both good and representable.
\end{corollary}

{\bf Proof.} Set $S=R/\Ann_RM$. Then $M$ possesses the structure
of an $S$-module in a natural way. A subset $N$ of $M$ is an
$R$-submodule of $M$ if and only if it is an $S$-submodule of $M$.
Thus it is straightforward to see that $M$ is good (resp.
representable) as an $R$-module if and only if it is good (resp.
representable) as an $S$-module. Now the conclusion follows by
1.1. $\Box$

\begin{proposition}  Let $R$ be a Noetherian ring. The following
statements are equivalent:\\
i) Every nonzero submodule of any representable $R$-module is
representable.\\
ii) $R$ is Artinian.
\end{proposition}

{\bf Proof.} $(ii)\Rightarrow (i)$ is clear by 1.1.\\
$(i)\Rightarrow (ii)$ By [{\bf 5}], any nonzero injective module
over a commutative Noetherian ring is representable. Since any
$R$-module can be embedded in an injective $R$-module, it follows
that all nonzero $R$-modules are representable. Therefore by the
implication $(iv)\Rightarrow (i)$ of 1.1, it follows that $R$ is
Artinian. $\Box$

A commutative ring $R$ is said to be Von Neumann regular, if for
each element $a\in R$, there exists $b\in R$ such that $a=a^2b$.
In [{\bf 1}, Theorem 2.3], it is shown that over a commutative Von
Neumann regular ring R every nonzero submodule of a representable
$R$-module is representable. Since commutative Artinian rings are
Noetherian, we can deduce the following result, by 2.4.

\begin{corollary}  Let $R$ be a commutative  Von Neumann regular
ring. Then $R$ is Noetherian if and only if it is Artinian.
\end{corollary}


\end{document}